\documentclass[letterpaper, 10 pt, conference]{ieeeconf}  

\IEEEoverridecommandlockouts                              

\usepackage{algorithmic}
\usepackage{amsmath, amssymb, amsxtra, mathrsfs, bm}
\usepackage{bbm}
\usepackage{color}
\usepackage{accents}
\usepackage{epsfig,here,graphicx}
\title{\LARGE \bf
Optimal Control of Transient Flows in Pipeline Networks with Heterogeneous Mixtures of Hydrogen and Natural Gas
}

\author{Luke S. Baker$^{1}$, Saif R. Kazi$^{2}$, Rodrigo B. Platte$^{1}$, and Anatoly Zlotnik$^{2}$
\thanks{*This study was supported by the U.S. Department of Energy's Advanced Grid Modeling (AGM) project ``Dynamical Modeling, Estimation, and Optimal Control of Electrical Grid-Natural Gas Transmission Systems'', as well as LANL Laboratory Directed R\&D project ``Efficient Multi-scale Modeling of Clean Hydrogen Blending in Large Natural Gas Pipelines to Reduce Carbon Emissions''. Research conducted at Los Alamos National Laboratory is done under the auspices of the National Nuclear Security Administration of the U.S. Department of Energy under Contract No. 89233218CNA000001.}
\thanks{$^{1}$Luke Baker and Rodrigo Platte are with the School of Mathematical and Statistical Sciences at Arizona State University,         Tempe, Arizona, 85281; {\tt \{lsbaker1,rplatte\}@asu.edu}.}
\thanks{$^{2}$Saif Kazi and Anatoly Zlotnik are in the Applied Mathematics \& Plasma Physics Group at Los Alamos National Laboratory, Los Alamos, New Mexico, 87545; {\tt \{skazi,azlotnik\}@lanl.gov}. }
}

\begin{document}

\maketitle
\thispagestyle{empty}
\pagestyle{empty}



\begin{abstract}
We formulate a control system model for the distributed flow of mixtures of highly heterogeneous gases through large-scale pipeline networks with time-varying injections of constituents, withdrawals, and control actions of compressors.  This study is motivated by the proposed blending of clean hydrogen into natural gas pipelines as an interim means to reducing end use carbon emissions while utilizing existing infrastructure for its planned lifetime.   We reformulate the partial differential equations for gas dynamics on pipelines and balance conditions at junctions using lumped elements to a sparse nonlinear differential algebraic equation system.  Our key advance is modeling the mixing of constituents in time throughout the network, which requires doubling the state space needed for a single gas and increases numerical ill-conditioning.  The reduced model is shown to be a consistent approximation of the original system, which we use as the dynamic constraints in a model-predictive optimal control problem for minimizing the energy expended by applying time-varying compressor operating profiles to guarantee time-varying delivery profiles subject to system pressure limits.  The optimal control problem is implemented after time discretization using a nonlinear program, with validation of the results done using a transient simulation.  We demonstrate the methodology for a small test network, and discuss scalability and potential applications.
\end{abstract}

\section{Introduction}

Transportation of natural gas through networks of large-scale transmission pipelines has been studied in steady-state \cite{wong1968optimization, percell1987steady, wu2000model, de2000gas, misra2014optimal} and transient \cite{abbaspour2008nonisothermal, osiadacz1984simulation, zlotnik2015optimal} operations with applications to the optimal control of compressor actuators.   In steady-state, the flow of gas in the network is balanced, so that inflows from processing plants and supply stations and outflows from withdrawal stations sum to zero.  Steady-state pipeline flows are described using simple time-invariant algebraic equations that relate pressure drop in the direction of flow to mass flow along each pipeline.  In the transient regime, computational complexity increases significantly because the flow in each pipeline cannot be modeled with simple algebraic equations but rather requires a system of nonlinear partial differential equations (PDEs) \cite{osiadacz2001comparison, thorley1987unsteady}.  Model reduction methods have been proposed to reduce the complexity of optimizing gas flows in networks \cite{rios2002reduction, grundel2013computing}.  Although natural gas is projected to be a primary fuel source through the year 2050 \cite{nalley2022annual}, worldwide economies have invested in transitioning from fossil fuels such as natural gas to more sustainable resources.  Hydrogen is a potential candidate, which, because it does not produce carbon dioxide when burned, is considered to have the potential to address climate change \cite{salvi2015sustainable}.  Natural gas pipeline operation and management protocols may be modified to transport mixtures of natural gas and hydrogen.  Recent studies indicate that natural gas pipelines can safely and effectively transport mixtures of up to 20\% hydrogen by volume \cite{gotz2016renewable, ozturk2021comprehensive}.  However, the complexity of modeling steady-state and transient flows, and thus designing and operating pipelines, is compounded with the injection of hydrogen \cite{van2004math, melaina2013blending}.

Natural gas and hydrogen have significantly different physical and chemical properties.  Hydrogen is less dense than natural gas, and the speed of sound through hydrogen is roughly four times as large as that of natural gas.  Viscosity, velocity, density, pressure, and energy of the gas mixture all vary with varying hydrogen concentration \cite{abd2021evaluation, blacharski2016effect}, and these directly affect the transportation of the mixture \cite{janusz2019hydrogen}.  Numerical simulations have been performed to demonstrate various effects on steady-state and transient-state flows of mixtures of hydrogen and natural gas in pipeline networks \cite{uilhoorn2009dynamic, chaczykowski2018gas, guandalini2017dynamic, elaoud2017numerical, hafsi2019computational, fan2021transient, agaie2017reduced}.  The method of characteristics was used in the numerical simulation of transient flows on cycle networks with homogeneous flow mixtures \cite{elaoud2017numerical}.   Another recent study investigates gas composition tracking using a moving grid method and an implicit backward difference method \cite{chaczykowski2018gas}.  It was shown that both methods of tracking perform well, but the implicit difference method may lose some finer detail in the response due to numerical diffusion.  A finite element method using COMSOL Multiphysics was also developed \cite{hafsi2019computational}.  That study considers the effects of hydrogen concentration on the compressibility factor of the mixture and its relation with pressure.  Moreover, there the authors demonstrate that pressure may exceed pipeline limitations in the transient evolution of flow and that the likelihood of this happening increases proportionally with increasing hydrogen concentration.

In contrast to pure natural gas, few studies have examined optimization of steady-state and transient operations of mixtures of hydrogen and natural gas in networks.  To our knowledge, there are no results on the optimal control of transient flows of heterogeneous mixtures of gases in pipelines or networks of pipelines.  Optimal control of compressor actuators for transport of pure natural gas typically seeks to minimize the cost of running compressors while being subjected to PDE flow dynamics, nodal pressure and nodal flow balance constraints, and inequality box constraints that limit the pressure throughout the network \cite{rachford2000optimizing}.  Other formulations may use an objective function that maximizes economic value \cite{zlotnik2019scheduling}.    When transients are sufficiently slow, a friction-dominated approximation may be made \cite{sundar2019tcst}, and this was shown to be valid in the regime of normal pipeline operations \cite{misra2020monotonicity}.  We use friction-dominated modeling to simplify the reduced modeling in the heterogeneous gas setting.



In this study, we formulate a control system model for transporting heterogeneous mixtures of gases through pipeline networks of general form, and extend optimal control problems for gas pipeline flow to this setting. Our key advance is modeling the mixing of constituents in time throughout the network, which requires doubling the state space needed for a single gas and increases numerical ill-conditioning.  This enables the formulation and solution of optimal control problems in which constituent gases may be injected at different points in the network at varying concentrations, e.g., the addition of 100\% hydrogen at certain nodes. An algorithm is implemented to obtain solutions, and the results are demonstrated on a small test network that includes a cycle.  

The remainder of the manuscript is organized as follows.  The governing equations for the flow of mixtures of gases in a network are presented in Section \ref{flow}.  In Section \ref{discretization}, an endpoint discretization method is employed to reduce the system of PDEs to a system of ordinary differential equations (ODEs).  There, we show that the discretization method is consistent and results in the equations for natural gas only in the case of zero hydrogen injection, and yields the steady-state equations when supply and withdrawal are held constant.  Section \ref{implementation} describes time-discretization of the optimal control problem that yields a nonlinear program (NLP).  The NLP is solved for a test network in Section \ref{study}, and we discuss applications in Section \ref{sec:conc}.

\section{Network Flow Control Formulation}  \label{flow}

We begin by defining notation that will be used in the study.  A gas network is modeled as a connected and directed graph $(\mathcal E, \mathcal V)$ that consists of edges (pipelines) $\mathcal E =\{1,\dots,E\}$ and nodes (junctions) $\mathcal V=\{1,2,\dots,V\}$, where $E$ and $V$ denote the cardinalities of the sets.  It is assumed that the nodes and edges are ordered within their sets according to their integer labels.  The symbol $k$ is reserved for identifying edges in $\mathcal E$ and the symbols $i$ and $j$ are reserved for identifying nodes in $\mathcal V$.   Supply nodes $ \mathcal V_s\subset \mathcal V$ and withdrawal nodes $ \mathcal V_w\subset \mathcal V$ are assumed to be disjoint sets that partition $\mathcal V$, i.e., $\mathcal V_s \cup \mathcal V_w = \mathcal V$ and $\mathcal V_s \cap \mathcal V_w = \emptyset$. It is assumed that supply nodes are ordered in $\mathcal V$ before withdrawal nodes so that $i<j$ for all $i\in \mathcal V_s$ and $j\in \mathcal V_w$.  The graph is directed by judiciously assigning a positive flow direction along each edge.  It is assumed that gas flows in the positive oriented direction of an edge so that the mass flux and velocity of the gas are positive quantities.  The notation $k:i\mapsto j$ means that edge $k\in \mathcal E$ is directed from node $i\in \mathcal V$ to node $j\in \mathcal V$.  For each node $j\in \mathcal V$, we define (potentially empty) incoming and outgoing sets of pipelines by $_{\mapsto} j=\{k\in \mathcal E| k:i\mapsto j\}$ and $j_{\mapsto}=\{k\in \mathcal E| k:j\mapsto i \}$, respectively.  

The transportation of the mixture of hydrogen and natural gas is modeled as a simplification of the isothermal Euler equations.  For each pipe $k\in \mathcal E$, the flow variables are natural gas density $\rho_k^{(1)}(t,x)$, hydrogen density $\rho_k^{(2)}(t,x)$, and mass flux $\varphi_k(t,x)$ of the mixture, with $t\in [0,T]$ and $x\in[0,\ell_k]$, where $T$ denotes the time horizon and $\ell_k$ denotes the length of the pipe. Assuming zero inclination and sufficiently slow transients, the flow through edge $k$ is governed by the friction-dominated PDEs
\begin{eqnarray}
\partial _t \rho_k^{(m)} +\partial_x \left(\frac{\rho_k^{(m)}}{\rho_k^{(1)}+\rho_k^{(2)}} \varphi_k \right) \!\!&=&\!\!0, \label{eq:flow1} \\
\partial_x\left(\sigma^2_1\rho_k^{(1)}+\sigma^2_2\rho_k^{(2)} \right) \!\!&=&\!\! -\frac{\lambda_k}{2D_k}\frac{ \varphi_k|  \varphi_k|}{\rho_k^{(1)}+\rho_k^{(2)}}, \label{eq:flow2} \quad
\end{eqnarray}
where \eqref{eq:flow1} is defined for each constituent $m=1$ and $m=2$.  Superscripts ``1" and ``2" attached to a gas variable are conserved for identifying natural gas and hydrogen variables, respectively.  The parameters for each pipeline $k\in \mathcal E$ are diameter $D_k$, friction factor $\lambda_k$, speed of sound through natural gas $\sigma_1$, and speed of sound through hydrogen gas $\sigma_2$.  In the above dynamic equations, the compressibility factors of the gasses are assumed to be constants so that the equations of states are ideally given by $p^{(m)}_k=\sigma^2_m\rho_k^{(m)}$, where $p_k^{(m)}$ is the partial pressure.  
The summation of partial pressures results in the equation of state $p_k=(\sigma^2_1\eta_k^{(1)}+\sigma^2_2\eta_k^{(2)})\rho_k$, where $p_k=(p^{(1)}_k+p^{(2)}_k)$ is the total pressure, $\rho_k=(\rho^{(1)}_k+\rho^{(2)}_k)$ is the total density, $\eta_k^{(1)}=\rho_k^{(1)}/\rho_k$ is the concentration of natural gas, and $\eta_k^{(2)}=\rho_k^{(2)}/\rho_k$ is the concentration of hydrogen.  The concentration as defined here refers to mass fraction, so that the volumetric fraction of hydrogen in the mixture is substantially greater.   

\begin{figure}[!t]
\centering
\includegraphics[width=.75\linewidth]{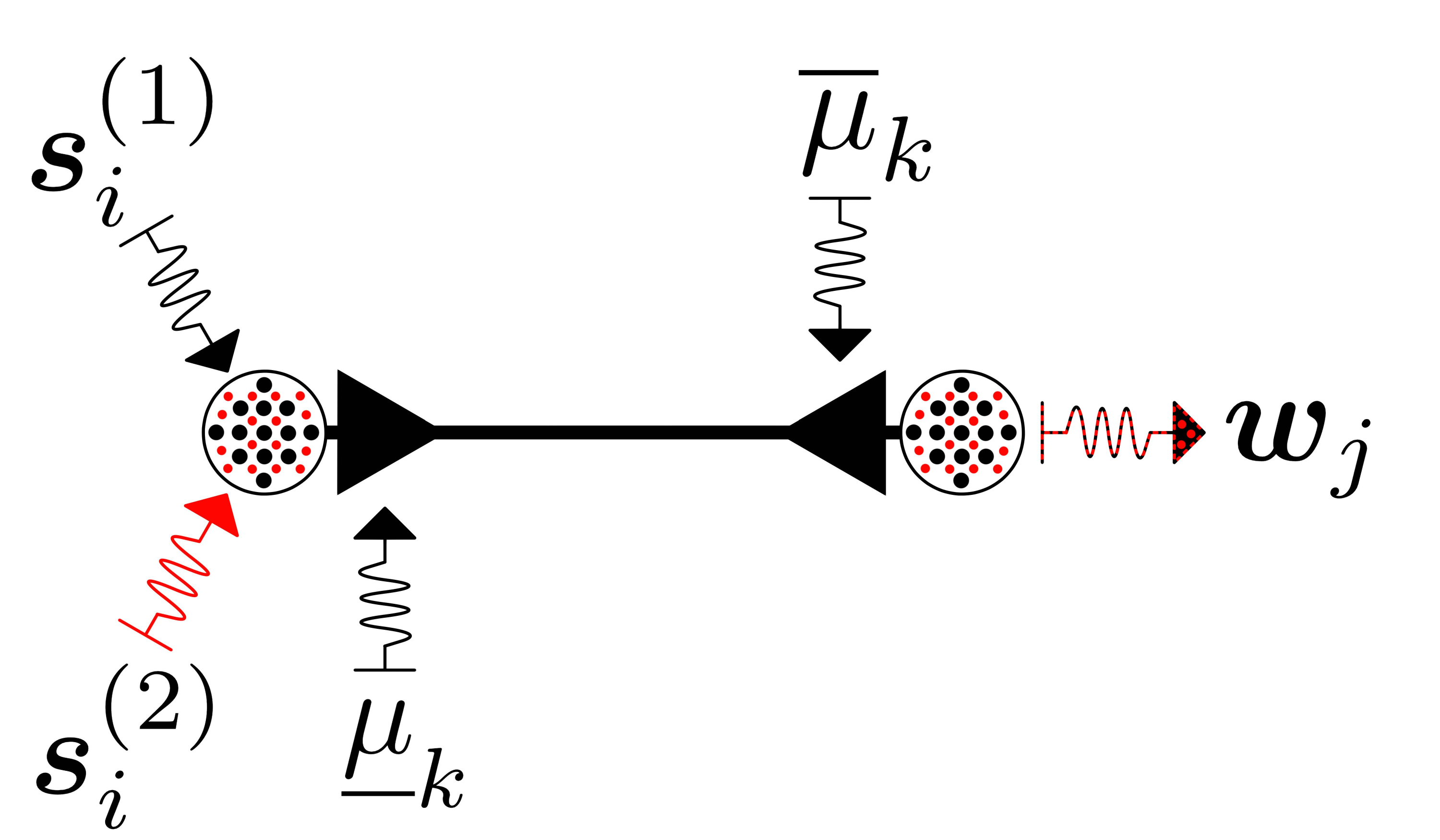}
\caption{Configuration of the pipeline segment $k:i\mapsto j$ with $i\in \mathcal V_s$ and $j\in \mathcal V_w$.}
\label{fig:pipe_conf}
\end{figure}

Friction forces between the interior wall of a pipe and gas flowing through it cause pressure to decrease in the direction of flow, as reflected in the momentum equation \eqref{eq:flow2}.  Compressor stations receive gas at low pressure and reduce its volume to increase its pressure to levels required for transportation and customer contracts. In addition to compressors, regulators are installed to reduce the pressure of the received gas to within limits that are compatible with lower pressure distribution systems.  For convenience, we assume that a compressor is located at the inlet and a regulator is located at the outlet of each pipeline with respect to the prescribed positive flow direction. For each pipeline $k\in \mathcal E$, compression and regulation are modeled with time-varying multiplicative compressor ratio $\underline\mu_k(t)\ge1$ and regulator ratio $\overline \mu_k(t)\ge1$, with orientations illustrated in Figure \ref{fig:pipe_conf}.

Natural gas and hydrogen are injected into the network at each supply node $i\in \mathcal V_s$ with specified time-varying profiles of natural gas density $ \bm s_i^{(1)}(t)$ and hydrogen density $\bm s_i^{(2)}(t)$.   Alternatively, pressure $(\bm p_s)_i=(\sigma_1^2 \bm s_i^{(1)}+\sigma_2^2 \bm s_i^{(2)})$ and concentration $ \bm \alpha_i^{(m)}=\bm s_i^{(m)}/(\bm s_i^{(1)}+\bm s_i^{(2)})$ may be specified at slack nodes $i\in \mathcal V_s$ instead of the constituent densities.  Gas is withdrawn downstream at each withdrawal node $j\in \mathcal V_w$ with specified time-varying mass flux $\bm w_j(t)$.  For $m=1,\,2$ and all $j\in \mathcal V_w$, define nodal density variables $\bm \rho_j^{(m)}(t)$ and nodal concentration variables $\bm \eta_j^{(m)}(t)$.  All of the nodal quantities in this study are identified with bold symbols.  Inlet and outlet edge variables are defined by attaching underlines below and overlines above the associated edge variables, respectively.  For example, $\underline \varphi_k(t)= \varphi_k(t,0)$ and $\overline \varphi_k(t)= \varphi_k(t,\ell_k)$.  The boundary conditions for the flow of the mixture are defined for $m=1,\,2$ by
\begin{eqnarray}
 \underline \rho_{k}^{(m)}&=&\underline\mu_{k} \bm s_i^{(m)},\qquad 
\overline \rho_{k}^{(m)}=\overline \mu_{k} \bm \rho_{j}^{(m)}, \label{eq:bc1} \\
\underline{\rho}_{k}^{(m)}&=&\underline\mu_{k} \bm \rho_{i}^{(m)},\qquad 
\overline \rho_{k}^{(m)}=\overline \mu_{k} \bm \rho_{j}^{(m)}, \label{eq:bc2} \\
\bm \eta_j^{(m)} \bm w_j&=&\sum_{k\in _{\mapsto}j}  \overline \eta_k^{(m)} \overline\varphi_{k} 
 -\sum_{k\in j_{\mapsto}} \underline \eta_k^{(m)} \underline\varphi_{k} ,  \label{eq:bc}
\end{eqnarray}
where  \eqref{eq:bc1} is defined for  $k:i\mapsto j$ with $i\in \mathcal V_s$, \eqref{eq:bc2} is defined for $k:i\mapsto j$ with $i\in \mathcal V_w$, and \eqref{eq:bc} is defined for $j\in \mathcal V_w$.   The configuration of the boundary conditions in a pipeline segment is depicted in Fig. \ref{fig:pipe_conf} (see Fig. \ref{fig:network} for a network example).  The conditions in \eqref{eq:bc1}-\eqref{eq:bc2} represent the effects of compression and regulation, and the conditions in \eqref{eq:bc} represent the conservation of mass flow of each constituent through withdrawal nodes.  It is assumed that the final operating state returns to its initial state, resulting in periodic temporal constraints
\begin{equation}
\rho_{k}^{(m)}(0,x)=\rho_{k}^{(m)}(T,x). \label{eq:ic}
\end{equation}
for all $k\in \mathcal E$ and $x\in [0,\ell_k]$.  Periodicity in time requires the parameters $\bm s_i^{(m)}(t)$, $\bm w_j(t)$, $\underline \mu_k(t)$, and $\overline \mu_k(t)$ to be periodic with period $T$.  We assume that the boundary conditions are smooth, slowly-varying, and bounded in their respective domains to ensure the existence of a smooth, slowly-varying, bounded solution.  The flow of the mixture of gases in the network is defined by the initial-boundary value system of PDEs \eqref{eq:flow1}-\eqref{eq:ic}.

Gas network operators require pressure, compression, and regulation to be within satisfactory limitations to ensure the safety of transportation and the quality of gas delivered to customers.  These limitations are modeled for all $k\in \mathcal E$ with inequality constraints of the form
\begin{equation}
 p^{\min}_k \le  \sigma_1^2\rho_{k}^{(1)}+\sigma_2^2\rho_{k}^{(2)} \le p^{\max}_k, \quad 1 \le  \underline\mu_{k}, \, \overline\mu_{k} \le 2, \label{eq:ineq}
\end{equation}
where $p^{\min}_k$ and $p^{\max}_k$ are specified bounds on pressure for each pipeline $k\in \mathcal E$. Compression $\underline\mu_{k}$ and regulation $\overline \mu_{k}$ are the control actuators in the network that are designed to minimize total consumption.  Because pressure down-regulation does not consume considerable energy, its proportion of total energy used for operating a pipeline may be neglected.  The total energy required for compression is given by
\begin{equation}
J = \sum_{k\in \mathcal E} \int_0^T c_{k}| \underline \varphi_{k}(t)|\left(  (\underline\mu_{k}(t))^{(\nu-1)/\nu}-1 \right)dt, \label{eq:J}
\end{equation}
where $c_{k}$ is related to the efficiency of the compressor $\underline \mu_{k}$ and $\nu$ is the isentropic exponent \cite{ivan2005} (which is assumed to be a weighted average of those of natural gas and hydrogen with weights equal to their respective mean concentration injections).   The continuous optimal control problem is 
\begin{equation} \label{eq:ocp1}
\begin{array}{ll}
\text{min}   \quad &J \triangleq \text{compressor energy in \eqref{eq:J}}, \\
\text{s.t.} \quad  & \text{dynamic constraints: \eqref{eq:flow1}-\eqref{eq:flow2}},  \\
 &\text{boundary conditions: \eqref{eq:bc1}-\eqref{eq:bc}},  \\
 &\text{temporal constraints: \eqref{eq:ic}},  \\
 &\text{inequality constraints: \eqref{eq:ineq}}. 
 \end{array}
\end{equation}
The decision variables are partial densities, mass fluxes, compressor ratios, and regulator ratios throughout the network.

\section{Network Flow Control Discretization} \label{discretization}

The intial-boundary value system of PDEs described in the previous section will be discretized in space to obtain an initial-value system of ODEs.   Spatial discretization is formalized by refining the graph of the gas network.  A graph refinement $(\hat{\mathcal E},\hat{\mathcal V})$ of the graph $(\mathcal E, \mathcal V)$ is made by adding auxiliary nodes to $\mathcal V$ that subdivide the edges of $\mathcal E$ so that $\ell_k\le \ell$ for all $k\in \hat{\mathcal{E}}$, where $\ell\le 10$ (km) is sufficiently small \cite{grundel2013computing}.  The refined graph inherits the prescribed direction of the parent graph.
For sufficiently fine network refinement, the relative difference of the flow variables between adjacent nodes is small in magnitude by continuity of the flow variables.  We assume that the graph has been sufficiently refined and that the hats may be omitted moving forward.

The system of ODEs is obtained by integrating the dynamic equations in \eqref{eq:flow1}-\eqref{eq:flow2} along the length of each refined pipeline segment so that
\begin{eqnarray*}
\int_0^{\ell}\partial _t \rho^{(m)}dx = -\int_0^{\ell}\partial_x \left(\frac{\rho^{(m)}}{\rho^{(1)}+\rho^{(2)}} \varphi \right)dx,  \\
\int_0^{\ell}\partial_x\left(\sigma^2_1\rho^{(1)}+\sigma^2_2\rho^{(2)} \right)dx
= -\frac{\lambda}{2D}\int_0^{\ell}\frac{ \varphi|  \varphi|}{\rho^{(1)}+\rho^{(2)}}dx,
\end{eqnarray*}
where edge subscripts have been removed for readability. The above integrals of space derivatives are evaluated using the fundamental theorem of calculus.  The remaining integrals are evaluated by approximating pipeline density with outlet density and pipeline flux with inlet flux.  These approximations are independent of $x$ and may be factored out of the integrals.  The above equations become
\begin{eqnarray}
\ell\dot {\overline\rho}^{(m)} &=&\underline \eta^{(m)} \underline\varphi-\overline \eta^{(m)} \overline\varphi, \label{dis1}
\\
\sum_{n=1}^{2}\sigma_n^2\left( \overline\rho^{(n)}-\underline \rho^{(n)} \right) &=& -\frac{\lambda \ell}{2D}\frac{\underline \varphi \left| \underline \varphi \right|}{\overline\rho^{(1)}+\overline \rho^{(2)}},  \label{dis2}
\end{eqnarray}
where a dot above a variable represents the time-derivative of the variable.  We now write the discretized system in matrix form.  Define the $E\times E$ diagonal matrices $L$ and $K$ with diagonal entries $L_{kk}=\ell_{k}$ and $K_{kk}=\lambda_k/(2D_k)$.  Define the time-varying (transposed) incidence matrix $M$ of size $E \times V$ componentwise by
\begin{eqnarray}
M_{ki}=  
\begin{cases}
\overline \mu_{k}(t), & \text{ edge $k\in _{\mapsto}i$ enters node $i$},
 \\
-\underline\mu_{k}(t), & \text{ edge $k\in i_{\mapsto}$ leaves node $i$,}
\\
0, & \text{else.}
\end{cases}
\end{eqnarray}
Define the $E\times r$ submatrix $M_s$ of $M$ by the removal of columns $i\in \mathcal V_w$, the $E\times (V-r)$ submatrix $M_w$ of $M$ by the removal of columns $i\in \mathcal V_s$, and the positive and negative parts of $M_w$ by $\overline M_w$ and $\underline M_w$ so that $M_w=(\overline M_w +\underline M_w)/2$ and $|M_w|=(\overline M_w -\underline M_w)/2$, where $r$ denotes the number of supply nodes and $|A|$ denotes the componentwise absolute value of a matrix $A$. Define the signed matrices $Q_w=$ sign$(M_w)$, $\overline Q_w=$ sign$(\overline M_w)$, $\underline Q_w=$ sign$(\underline M_w)$, and similarly for $M_s$.  These signed matrices are well-defined by the inequalities in \eqref{eq:ineq}. 

Define inlet and outlet edge mass flux vectors by $\underline \varphi =(\underline \varphi_1,\dots,\underline \varphi_E)^T$ and $\overline \varphi =(\overline \varphi_1,\dots,\overline \varphi_E)^T$, and similarly for inlet and outlet edge concentrations.  Moreover, define the vectors $\bm \rho^{(m)}=(\bm \rho_{r+1}^{(m)},\dots, \bm \rho_{V}^{(m)})^T$, $\bm \eta^{(m)}=(\bm \eta_{r+1}^{(m)},\dots, \bm \eta_{V}^{(m)})^T$, and $\bm \alpha^{(m)}=(\bm \alpha_1^{(m)},\dots, \bm \alpha_r^{(m)})^T$, where the subscripts of the entries are indexed according to the node labels in $\mathcal V$.  Applying the above matrix definitions, the discretized equations in \eqref{dis1}-\eqref{dis2} together with the boundary conditions in \eqref{eq:bc1}-\eqref{eq:bc} become
\begin{eqnarray}
L\overline M_w \dot {\bm \rho}^{(m)} = \underline \eta^{(m)} \odot \underline \varphi -\overline \eta^{(m)}  \odot \overline \varphi , \label{Dis1} 
\\
\sum_{m=1}^2 \sigma_m^2 \left( M_w \bm \rho^{(m)}+M_s \bm s^{(m)} \right) = -\frac{LK (\underline \varphi \odot |\underline \varphi|)}{\overline M_w(\bm \rho^{(1)}+\bm\rho^{(2)})},  \label{Dis2} \\
\bm \eta^{(m)}\odot \bm w=\overline Q_w^T\left( \overline \eta^{(m)}  \odot \overline \varphi \right) +\underline Q_w^T\left( \underline \eta^{(m)} \odot \underline \varphi \right), \label{Dis3} 
\end{eqnarray} 
where $\odot$ is the Hadamard product, and the ratio of vectors on the right-hand-side of \eqref{Dis2} is understood to be componentwise.  It is assumed that regulators vary slowly so that the time derivative of $\overline M_w$ is insignificant, justifying its removal from \eqref{Dis1}.  Multiplying both sides of \eqref{Dis1} on the left by $\overline Q_w^T$ and using \eqref{Dis3}, we may combine \eqref{Dis1} and \eqref{Dis3} to form the equation $\overline Q_w^T L\overline M_w\dot {\bm \rho}^{(m)}=[Q_w^T(\underline \eta^{(m)} \odot \underline \varphi)- \bm \eta^{(m)} \odot \bm w]$, where we have used $Q_w=(\underline Q_w+\overline Q_w)$.  By the definitions of supply and withdrawal concentrations, the above equations become 
\begin{eqnarray}
\overline Q_w^T L\overline M_w\dot {\bm \rho}^{(m)}=Q_w^T[(|\underline Q_w| \bm \eta^{(m)}+|\underline Q_s| \bm \alpha^{(m)} )\odot \underline \varphi]
\nonumber \\
- \bm \eta^{(m)} \odot \bm w, \label{partial1} \\
\sum_{m=1}^2 \sigma_m^2 \left( M_w \bm \rho^{(m)}+M_s \bm s^{(m)} \right) = -\frac{LK ( \underline \varphi \odot |\underline \varphi|)}{\overline M_w(\bm \rho^{(1)}+\bm\rho^{(2)})}. \label{partial2}
\end{eqnarray} 
Periodic temporal constraints in \eqref{eq:ic} reduce to
\begin{eqnarray}
\bm \rho^{(m)}(0)=\bm \rho^{(m)}(T). \label{ic2}
\end{eqnarray}
Pressure, compression, and regulation inequality constraints in \eqref{eq:ineq} reduce to
\begin{equation}
 \bm p^{\min}_j \le  \sigma_1^2\bm \rho_{j}^{(1)}+\sigma_2^2\bm \rho_{j}^{(2)} \le \bm p^{\max}_j, \quad 1 \le  \underline\mu_{k}, \, \overline\mu_{k} \le 2, \label{Disineq}
\end{equation}
where $ \bm p^{\min}_j$ and  $\bm p^{\max}_j$ are specified bounds for each node $j\in \mathcal V_w$.  The reduced-model optimal control problem is formulated as
\begin{equation} \label{eq:Discop}
\begin{array}{ll}
\text{min}   \quad &J\triangleq \text{compressor energy in \eqref{eq:J}}, \\
\text{s.t.} \quad  & \text{dynamic constraints: \eqref{partial1}-\eqref{partial2}},  \\
  &\text{temporal constraints: \eqref{ic2}},  \\
 &\text{inequality constraints: \eqref{Disineq}}. 
\end{array} 
\end{equation}

We now present a few results on the discretization method.  Proposition 1 below shows that the discretized system in \eqref{partial1}-\eqref{partial2} approaches the continuous system in \eqref{eq:flow1}-\eqref{eq:flow2} in a single pipeline as the distance between adjacent nodes of the refined pipeline approaches zero.  Proposition 2 shows that the number of density variables in \eqref{partial1}-\eqref{partial2} reduces to half this number for homogeneous mixtures. Moreover, there we show that the discretized system reduces to the steady-state equations in the time-invariant setting.

{\bf Proposition 1.}  Consider a single pipeline of length $\ell$, and refine its graph as a chain connection of $E$ segments of uniform length $\Delta \ell=\ell/E$, diameter $D$, and friction factor $\lambda$.  Suppose the gas mixture is supplied to the pipeline at the inlet $\mathcal V_s=\{1\}$ with boundary conditions as in \eqref{eq:bc1} and withdrawn from only the outlet so that $\bm w_j=0$ for $j\not= E+1$.  Suppose for simplicity that there are no compressors or regulators.  Then the resulting system in \eqref{partial1}-\eqref{partial2} is a consistent spatial discretization of \eqref{eq:flow1}-\eqref{eq:bc}.

{\bf Proof.} The matrix $\overline M_w$ is the $E\times E$ identity matrix, $\underline M_w$ is the $E\times E$ lower off-diagonal matrix with nonzero entries $(\underline M_w)_{n+1,n}=-1$, and $\underline Q_s$ is an $E\times 1$ unit vector with one nonzero entry given by $(\underline Q_s)_{1}=-1$.  For the intermediate segment $n:n\mapsto (n+1)$ with $2\le n \le E-1$, the associated dynamics in \eqref{partial1}-\eqref{partial2} are given by
\begin{IEEEeqnarray*}{rcl}
 \dot {\bm \rho}^{(m)}_{n+1}+\frac{1}{\Delta \ell}\left(\bm \eta_{n+1}^{(m)}\underline \varphi_{n+1} -\bm \eta_{n}^{(m)}\underline \varphi_n\right) &=& \; 0 \\
\frac{1}{\Delta \ell}\sum_{m=1}^{2}\sigma_m^2\left( \bm \rho_{n+1}^{(m)}-  \bm \rho_{n}^{(m)}\right) &=& -\frac{\lambda}{2D}\frac{\underline \varphi_n |\underline \varphi_n|}{\bm \rho_{n+1}^{(1)}+\bm \rho_{n+1}^{(2)}}.
\end{IEEEeqnarray*}
Taking the limit $\Delta \ell \rightarrow 0$, the above equations approach the dynamics in \eqref{eq:flow1}-\eqref{eq:flow2}.  Similarly, as $\Delta \ell\rightarrow 0$, it can be shown that the first and last segments of the pipe reduce to the dynamics \eqref{eq:flow1}-\eqref{eq:flow2} with boundary conditions \eqref{eq:bc1}-\eqref{eq:bc}. $\square$

{\bf Proposition 2.}  The network system in \eqref{partial1}-\eqref{partial2} reduces to a system with total density as the only density variable if the concentration is homogeneous.  If the concentration of hydrogen is zero, then the system reduces further to the pure natural gas equations.  Furthermore, if the concentration, supply pressure, and withdrawal flux are time-invariant, yielding a time-invariant solution, then the system in \eqref{partial1}-\eqref{partial2} reduces to the steady-state balance laws.

{\bf Proof.}  Suppose that $\bm \alpha^{(2)}:=\bm \alpha$ and $\bm \eta^{(2)}:=\bm \eta$ are constant vectors that are known.  The equation of state may be written as $(\sigma_1^2\bm \rho^{(1)}+\sigma_2^2\bm \rho^{(2)})=\bm a^2 \odot \bm \rho$, where $\bm \rho=(\bm \rho^{(1)}+\bm \rho^{(2)})$ is the total density vector and $\bm a=(\sigma_1^2\bm (1-\bm \eta)+\sigma_2^2\bm \eta)^{1/2}$ is a generalized sound speed vector (the square-root operation is componentwise). Moreover, $(\sigma_1^2\bm s^{(1)}+\sigma_2^2\bm s^{(2)})=\bm b^2 \odot \bm s$, where $\bm s=(\bm s^{(1)}+\bm s^{(2)})$ and $\bm b=(\sigma_1^2\bm (1-\bm \alpha)+\sigma_2^2\bm \alpha)^{1/2}$. Because $\bm a$ and $\bm b$ are constant vectors, the superposition of \eqref{partial1}, for $m=1,\,2$, results in
\begin{align}
\overline Q_w^T L\overline M_w\dot {\bm \rho}&= Q_w^T\underline \varphi -  \bm w,  \label{reduced1} \\ \nonumber \\
M_w \left( \bm a^2 \odot \bm \rho\right)+M_s \left( \bm b^2 \odot \bm s \right) &= -\frac{LK ( \underline \varphi \odot |\underline \varphi|)}{\overline M_w\bm \rho}. \quad \label{reduced2}
\end{align} 
Partial densities may be determined using $\bm \rho^{(m)}=\bm \eta^{(m)}\odot \bm \rho$.  This shows that, under the assumption of constant concentration, the reduced system in \eqref{reduced1}-\eqref{reduced2} may be used to determine the solution to \eqref{partial1}-\eqref{partial2}.  If $\bm \alpha=0$ and $\bm \eta=0$, then the above equations reduce to the single gas endpoint discretization method with $\bm a=\bm b=\sigma_1$ \cite{himpe2021model}.  If $\bm w$, $\bm \eta$, $\bm s$, and $\bm \rho$ are constant, then the system \eqref{reduced1}-\eqref{reduced2} reduces to the Weymouth equations for a mixture of gases \cite{misra2014optimal}.
$\square$


\section{Implementation} \label{implementation}

The optimal control problem in \eqref{eq:Discop} may be expressed as
\begin{subequations}
\begin{align}
\min   \qquad         &\int_0^T \mathcal{F}(x(t),u(t))dt \\
\text{s.t.}  \qquad   &F \frac{d}{dt} [R^{(m)}x (t)]=f^{(m)} (x(t),u(t),p(t)), \\  
    & e(x(T),x(0),u(T),u(0))=0, \\
    & q(x(t),u(t)) \ge 0,
\end{align}
\end{subequations}
where $x=(\bm \rho^{(1)},\bm \rho^{(2)},\underline \varphi)^T$ is the state, $u=(\{\underline \mu_k, \overline \mu_k\})^T$ is the control, and $p= (\bm s^{(1)},\bm s^{(2)},\bm w )^T$ is a fixed (potentially time-varying) vector of parameters.  The function $e$ represents temporal periodic constraints, $q$ represents inequality constraints, $R^{(m)}$ is a row selector matrix that maps $x$ into $\bm \rho^{(m)}$, and $F=\overline Q_w^T L\overline M_w$.  A nonlinear program is obtained by discretizing the time interval $[0,T)$ into $N$ subintervals with equally-spaced collocation points $t_n=(n-1)T/N$ for $n=1,\dots,N$. 

The vector-valued functions $x(t)$, $u(t)$, and $p(t)$ are interpolated with piecewise-linear vector-valued functions $\hat{x}(t)$, $\hat u(t)$, and $\hat p(t)$, respectively.  For $n=1,\dots,N$, the function $\hat x(t)$ is defined for $t\in [t_n,t_{n+1})$ by
\begin{equation*}
    \hat x(t)=x(t_n)+\frac{x(t_{n+1})-x(t_n)}{T/N}(t-t_n),
\end{equation*}
where $x(t_{N+1})=x(t_1)$ by the assumption of periodicity.  The functions $\hat u(t)$ and $\hat p(t)$ are defined similarly.   The integral in the objective function is approximated using the trapezoidal rule, resulting in
\begin{equation*}
    \int_0^T \mathcal{F}(x(t),u(t))dt \approx \sum_{n=1}^{N} \frac{T}{N}\mathcal{F}(\hat x(t_n),\hat u(t_n)).
\end{equation*}
The time derivative of $R^{(m)}x(t)$ is approximated using the finite difference $d/dt[R^{(m)} x(t)] \approx N/T( R^{(m)}x(t_{n+1})- R^{(m)}x(t_{n}))$.  Define the $N\times N$ differentiation matrix $D$ componentwise by $D_{n,n}=-T/N$, $D_{n,n+1}=T/N$, and $D_{N,1}=T/N$.  The differentiation matrix includes a periodic constraint to the vector on which it operates.  Define the stacked state vector $X=(\hat x(t_1),\dots,\hat x(t_N))^T$, input $U=(\hat u(t_1),\dots,\hat u(t_N))^T$, parameter $P=(\hat p(t_1),\dots,\hat p(t_N))^T$, and extend $D$, $F$, $R^{(m)}$, and all of the other matrices using Kronecker products with identity matrices of appropriate dimensions. Then the nonlinear program may be written as
\begin{subequations} \label{NLP}
\begin{align}
\min   \qquad         &\sum_{n=1}^{N} \frac{T}{N}\mathcal{F}(X_nX, U_nU) \\
\text{s.t.}  \qquad   &F DR^{(m)}X_nX=f^{(m)} (X_nX,U_nU,P_nP), \nonumber \\
                       &m=1,2,\; n=1,\dots,N, \\  
    & q(X_nX,U_nU) \ge 0, \; n=1,\dots,N,
\end{align}
\end{subequations}
where $X_n$, $U_n$, and $P_n$ are row selector matrices that satisfy $X_nX=\hat x(t_n)$, $U_nU=\hat u(t_n)$, and $P_nP=\hat p(t_n)$ for $n=1,\dots,N$.  

Solution of the NLP in \eqref{NLP} is implemented in Matlab with the interior-point algorithm using the function \verb"fmincon", and is evaluated on a MacBook Air 8‑core CPU with 8GB of unified memory. The gradient of the objective and Jacobian of the constraints are supplied to the function for improved performance.  The Hessian of the Lagrangian function is set to the default finite-difference approximation.  Optimal control of compression and regulation is obtained from the optimal solution $U=U^*$.  The optimized time-series for compressors and regulators, and any specified parameters, are linearly interpolated to provide control functions to the ODEs \eqref{partial1}-\eqref{partial2}.  The system is simulated in Matlab using the function \verb"ode15s" for validation of the solution and an improved prediction of pressure and mass flux.  The steady-state solution is used as the starting point for optimization, and the initial state of the optimal solution is used as the initial condition for simulation.  In the following, we distinguish between the solution of the optimization problem \eqref{NLP} and the solution of the ODEs \eqref{partial1}-\eqref{partial2} that are driven by optimal compression and regulation. We consider this comparison an important validation of the presented optimal control scheme, because feasibility of the coarsly discretized physical system in the optimal solution of the NLP \eqref{NLP} does not necessarily guarantee that the control solution obtained by solving \eqref{NLP} results in that same physical solution in a simulation with controlled error.  For the case studies in the next section, the two solutions are compared using the average $L^2$ norm of the relative difference given by
\begin{equation}
    \frac{1}{\hat{E}} \sum_{k\in \hat{ \mathcal E}} \left( \frac{1}{T}\int_0^T  \left( 2\frac{\underline \varphi_k(t) -\underline \phi_k(t)}{\underline \varphi_k(t) +\underline \phi_k(t)} \right)^2 dt \right)^{1/2}\times 100, \label{2_norm}
\end{equation}
where $\underline \varphi_k$ is the optimized flux in edge $k\in \hat{\mathcal E}$, $\underline \phi_k$ is the simulated flux, and $\hat{ E}$ is the cardinality of $\hat{\mathcal E}$.  In addition, the maximum absolute relative difference is also documented as
\begin{equation}
    \max_{k\in \hat{ \mathcal E}} \left( \max_{t\in [0,T]} \left| 2\frac{\underline \varphi_k(t) -\underline \phi_k(t)}{\underline \varphi_k(t) +\underline \phi_k(t)}\right| \right)\times 100. \label{max_norm}
\end{equation}
Similar metrics are used for the difference between optimized and simulated pressure trajectories.

\begin{figure}[!t]
\centering
\includegraphics[width=\linewidth]{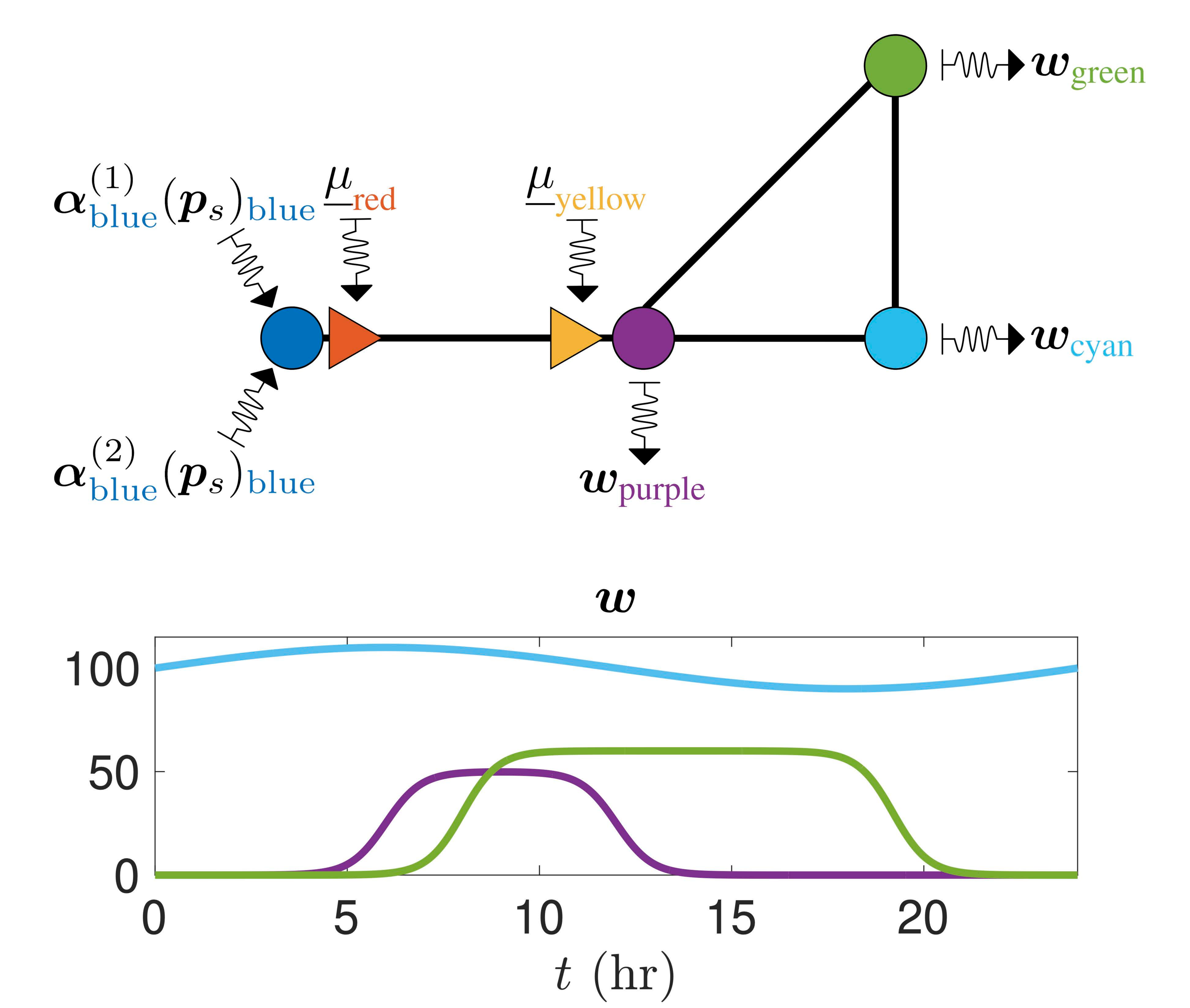}
\caption{Top: Configuration of the network.  Pipeline dimensions: blue to purple = 50 km, purple to cyan = 30 km, purple to green = 20 km, and green to cyan is 30 km.  Bottom: Withdrawal flow profiles, color-coded to correspond to associated nodes.  The network has two compressors, indicated by triangles, located at the start and end of the pipe directed from blue to purple.}
\vspace{1ex}
\label{fig:network}
\end{figure}

\section{Case Study} \label{study}

The optimal control algorithm is demonstrated on a cyclic network whose configuration and dimensions are shown on the top of Fig. \ref{fig:network}.  We create a refined network with a uniform discretization length of $\ell_k=10$ km for all $k\in \hat{\mathcal E}$.  The diameters and friction factors of the refined pipelines are uniform and equal to $D_k=0.5$ m and $\lambda_k=0.011$ for all $k\in \hat{\mathcal E}$.  The speeds of sounds of the gases are $\sigma_1=338.38$ m/s and $\sigma_2=4\sigma_1$.   We use $N=20$ time steps with $\nu=1.28$ and compressor efficiency values $c_1=c_5=D_k^4/T$ in \eqref{eq:J}. The minimum and maximum pressures in \eqref{Disineq} are $\bm p^{\min}_j=5$ MPa and $\bm p^{\max}_j=12$ MPa for all $j\in \hat{\mathcal V}_w$.  Discretization results in 780 optimization variables, 740 equality constraints, and 520 inequality constraints in the NLP \eqref{NLP}.

\begin{figure}[!t]
\centering
\hspace{-3ex}\includegraphics[width=\linewidth]{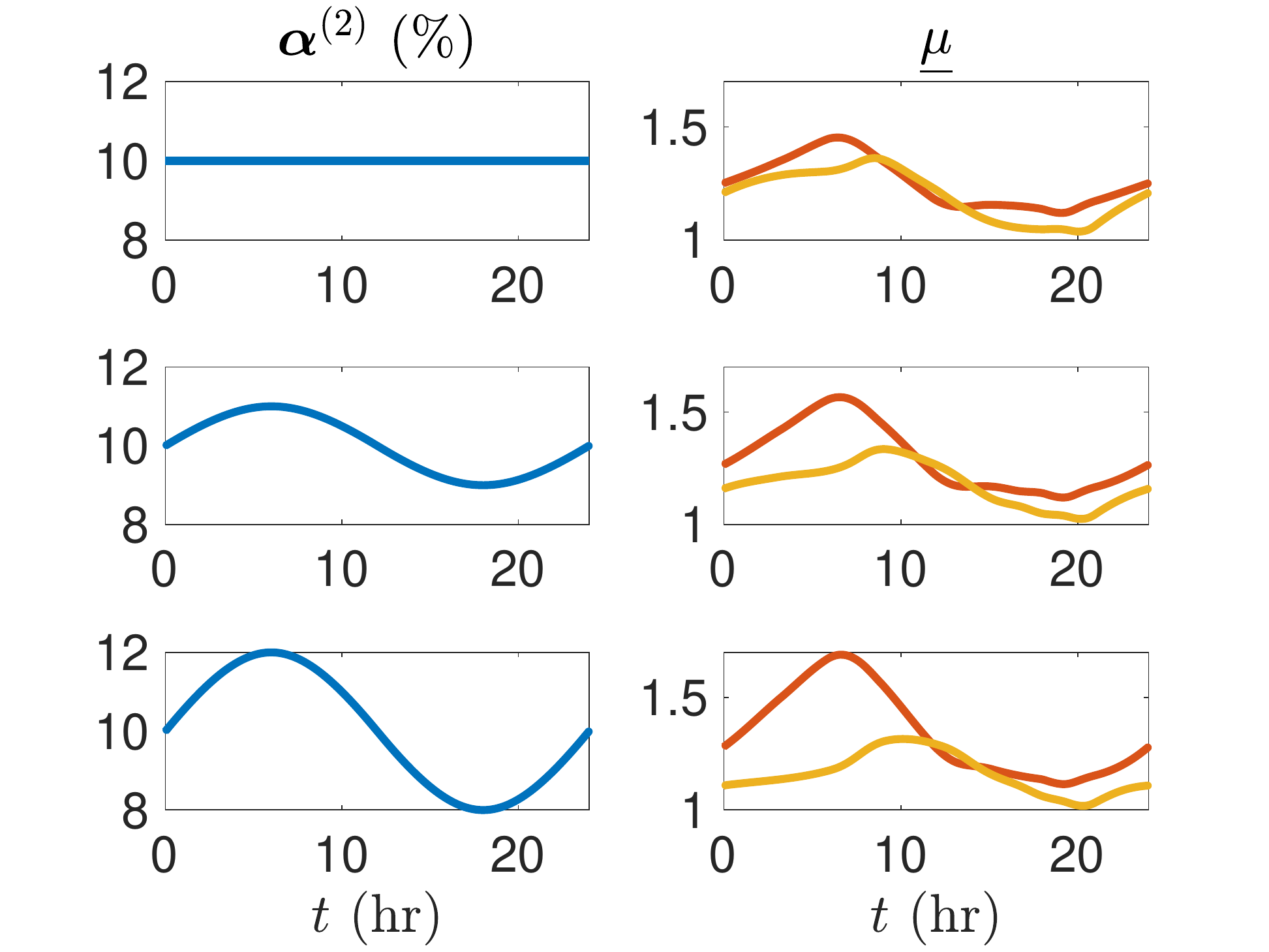}
\caption{Left column: Hydrogen concentration profiles at the supply node. Right column: Corresponding optimal compressor responses.}
\vspace{3ex}
\label{fig:concentrations}
\end{figure}

\begin{figure}[t!]
\centering
\hspace{-3ex}\includegraphics[width=\linewidth]{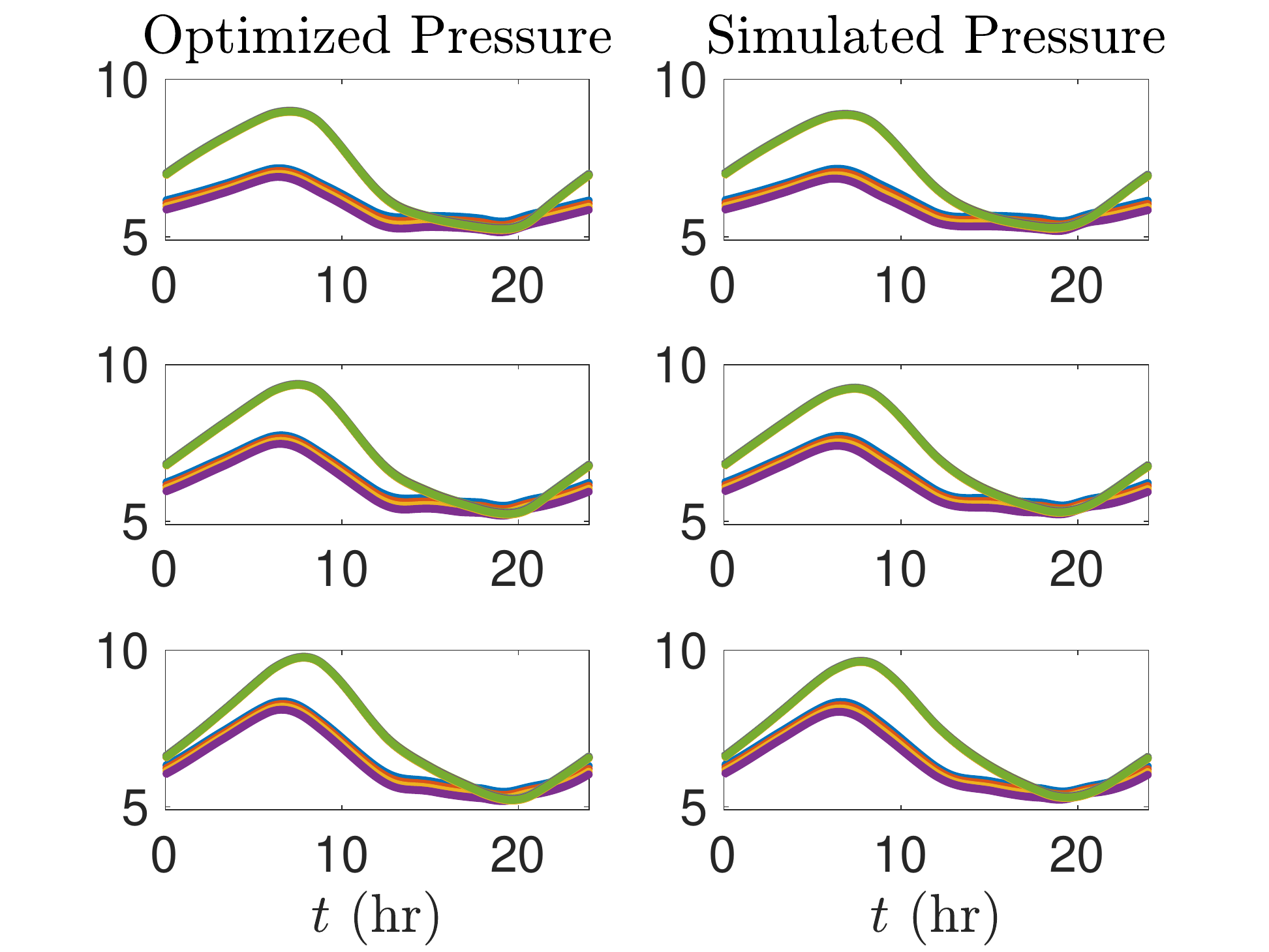}
\caption{Left column: Optimized pressure at refined withdrawal nodes driven by the respective concentration profiles from Fig. \ref{fig:concentrations}. Right column: Simulated pressure at refined withdrawal nodes driven by inlet concentration and optimal compression profiles from Fig. \ref{fig:concentrations}.  The average $L^2$ norm of the relative difference in pressures, using the metric \eqref{2_norm}, are approximately 0.769\%, 0.770\%, and 0.769\% for the top, middle, and bottom rows.  The maximum relative difference in pressures, using the metric \eqref{max_norm}, are approximately 2.154\%, 2.038\%, and 1.971\% for the top, middle, and bottom rows.  The combined computational times for optimization and simulation are approximately 10.0 s, 8.0 s, and 7.7 s for the top, middle, and bottom rows, respectively.}
\label{fig:pressures}
\end{figure}

\begin{figure}[t!]
\centering
\hspace{-3ex}\includegraphics[width=0.5\textwidth]{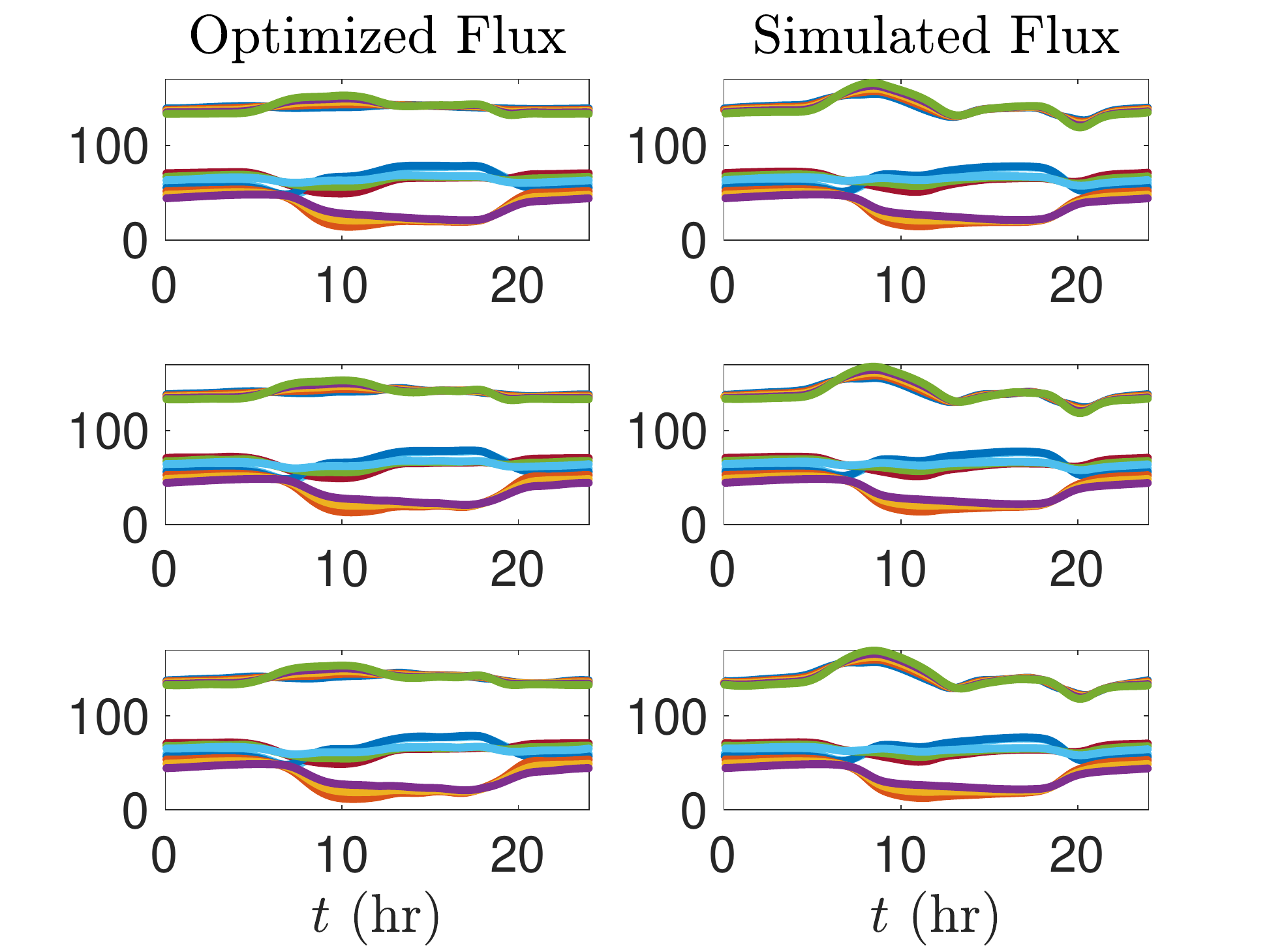}
\caption{Left column: Optimized mass flux in refined edges driven by the respective concentration profiles from Fig. \ref{fig:concentrations}. Right column: Simulated mass flux in refined edges driven by inlet concentration and optimal compression profiles from Fig. \ref{fig:concentrations}.  The $L^2$ relative difference metric values using \eqref{2_norm} are approximately 3.994\%, 4.608\%, and  5.258\% for the top, middle, and bottom rows.  The maximum relative difference values using metric \eqref{max_norm} are approximately 12.967\%, 16.713\%, and 21.509\%, respectively.}
\vspace{1ex}
\label{fig:massflux}
\end{figure}

The purple, green, and cyan nodes in the network graph in Figure \ref{fig:network} represent stations where gas is withdrawn with color-coordinated flow profiles depicted on the bottom of Fig. \ref{fig:network}.  The red and yellow triangles represent two compressor stations whose time-dependent operations are optimized in a  model-predictive manner.  The blue node is the supply station for a mixture of natural gas and hydrogen with a fixed pressure of $(\bm p_s)_{\text{blue}}=5$ MPa that is immediately boosted by the red compressor station.  We demonstrate three solutions for this network, each of which are subject to the same above boundary conditions but differ in the injected concentration of hydrogen at the supply node. The left column of Fig. \ref{fig:concentrations} depicts the specified hydrogen concentration profiles at the supply node for the three solutions and the right column shows the associated results for optimal compression ratios of the two color-coordinated compressor stations.  The total compressor energy values in as defined in equation \eqref{eq:J} are $J=0.787$, $J=0.824$, and $J=0.860$ (non-dimensionalized units) corresponding to the optimal compression ratios given in Fig. \ref{fig:concentrations} from top to bottom, respectively.  Figures \ref{fig:pressures} and \ref{fig:massflux} depict pressure and mass flux solutions, respectively, where the left-hand-side columns show the optimized solutions and the right-hand-side columns show the validating simulation.  The results demonstrate that minor variations in hydrogen concentration may have substantial effects on pressure and compressor activity.  In particular, the pressure trajectories in Figure \ref{fig:pressures} increase by approximately 2 MPa from $t=0$ to $t=8$ hours for a fixed 10\% hydrogen injection, while the change is approximately 3.5 MPa when the hydrogen injection concentration slowly varies from 8\% to 12\% over the 24 hour time horizon.  Observe that because this concentration is defined as the mass fraction, the volume fraction here is actually over 30\%, and exhibits variations of up to 10\%.  A detailed analysis of how hydrogen blending impacts energy transport capacity of gas pipeline networks is outside the scope of this study, in which we focus on modeling flow dynamics and demonstrating our optimal control approach.

\section{Conclusions} \label{sec:conc}

We synthesized a control system model for the distributed flow of mixtures of two gases with different physical properties through large-scale pipeline networks with time-varying injections, withdrawals, and control actions of compressors.  The motivation is to develop analysis methods to evaluate recent proposals for blending of clean hydrogen into natural gas pipelines as an interim means for carbon emissions reduction  that allows utilization of existing infrastructure for its planned lifetime \cite{raju2022}.   The partial differential equations for gas dynamics on pipelines and balance conditions at junctions are approximated using lumped elements to a sparse nonlinear differential algebraic equation system.  In contrast to previously published results, we model the mixing of constituents in time throughout the network.  The optimal control technique we develop is able to rapidly produce validated solutions, even though representing dynamics of the gas mixture requires doubling the state space with respect to models for a single gas, and worsens numerical conditioning.  We show that the reduced model is a consistent approximation of the original system, use it as the dynamic constraints in a model-predictive optimal control method for minimizing the energy expended by applying time-varying compressor operating profiles to guarantee time-varying delivery profiles subject to system pressure limits.  The optimal control problem is implemented after time discretization using a nonlinear program, with validation of the results done using a transient simulation.  

The developed control system model and computational optimal control scheme can be used to solve a variety of problem formulations for gas transport networks.  The objective function could be modified to reflect the economic value of pipeline transport, in terms of natural gas and hydrogen flow provided by suppliers, and energy received by consumers.  Including a price of carbon emissions mitigation due to replacement of natural gas with hydrogen could indicate optimal locations for integrating hydrogen supplies.

\bibliographystyle{IEEEtran}
\bibliography{sample.bib}

\end{document}